\pdfoutput=1
\documentclass[nohyperref]{article}

\usepackage{microtype}
\usepackage{graphicx}
\usepackage{subfigure}
\usepackage{booktabs} 

\usepackage{hyperref}

 \usepackage[accepted]{icml2022}

\usepackage{amsmath}
\usepackage{amssymb}
\usepackage{mathtools}
\usepackage{amsthm}
\usepackage{dsfont}

\usepackage[capitalize,noabbrev]{cleveref}

\theoremstyle{plain}
\newtheorem{theorem}{Theorem}[section]

\theoremstyle{definition}

\newtheorem{assumption}[theorem]{Assumption}
\theoremstyle{remark}

\usepackage[textsize=tiny]{todonotes}

\DeclareMathOperator*{\argmin}{arg\,min}

\icmltitlerunning{Two-Timescale Stochastic Approximation for Bilevel Optimisation Problems in Continuous-Time Models}

\begin{document}

\onecolumn
\icmltitle{Two-Timescale Stochastic Approximation for Bilevel \\ Optimisation Problems in Continuous-Time Models}



\icmlsetsymbol{equal}{*}

\begin{icmlauthorlist}
\icmlauthor{Louis Sharrock}{yyy}
\end{icmlauthorlist}

\icmlaffiliation{yyy}{Department of Mathematics, University of Bristol, Bristol, UK}

\icmlcorrespondingauthor{Louis Sharrock}{louis.sharrock@bristol.ac.uk}

\icmlkeywords{Stochastic Gradient Descent, Stochastic Approximation, Convergence Rate, Central Limit Theorem}

\vskip 0.3in



\printAffiliationsAndNotice{}  

\begin{abstract}
We analyse the asymptotic properties of a continuous-time, two-timescale stochastic approximation algorithm designed for stochastic bilevel optimisation problems in continuous-time models. 
We obtain the weak convergence rate of this algorithm 
in the form of a central limit theorem. We also demonstrate how this algorithm can be applied to several continuous-time bilevel optimisation problems. 
\end{abstract}

\section{Introduction}
In recent years, bilevel optimisation problems have received significant attention, owing no doubt in part to their relevance to a wide range of machine learning applications. These include meta-learning \cite{Franceschi2018,Bertinetto2019}, 
hyperparameter optimisation \cite{Feurer2019,Shaban2019}, 
and reinforcement learning \cite{Konda2003a,Khodadadian2021}. 
In this paper, we consider bilevel optimisation problems of the form
\begin{equation}
\min_{x\in\mathbb{R}^{d_1}} \Phi(x) = f(x,y^{*}(x))~~~\text{subject to}~~~y^{*}(x) \in \argmin_{y\in\mathbb{R}^{d_2}} g(x,y), \label{eq:bilevel}
\end{equation}
where $d_1,d_2\in\mathbb{N}$ are integers, and $f,g:\mathbb{R}^{d_1}\times\mathbb{R}^{d_2}\rightarrow\mathbb{R}$ are continuously differentiable functions. We will refer to $f$ as the \emph{upper-level} or \emph{outer-level} objective function, and $g$ as the \emph{lower-level} or \emph{inner-level} objective function. 

To solve this problem, 
we consider an 
 approach, first proposed by \cite{Hong2020}, in which approximate solutions to the outer- and inner-problems are updated simultaneously, based on stochastic estimates of the gradients of the two objective functions. This method belongs to broad class of stochastic optimisation algorithms known as \emph{two-timescale stochastic approximation} (TTSA), which consist of two coupled recursions evolving on different timescales \cite{Borkar1997}.
 This approach has found success in a wide variety of applications \cite{Heusel2017,Yang2019,Wu2020}.

The asymptotic properties of this approach, under various assumptions, have been well studied \cite{Tadic2003,Tadic2004}.
 More recently, its non-asymptotic properties have also received significant attention \cite{Dalal2018,Gupta2019,Kaledin2020,Doan2021a}. 
 While existing works provide an excellent insight, 
 they only analyse two-timescale algorithms in discrete time. In this article, we instead consider TTSA in continuous time, motivated by bilevel optimisation problems arising in continuous-time models  \cite{Doya2000,Pham2009,Sirignano2017a,Yildiz2021}. 

Even in the single timescale case, results on continuous-time stochastic approximation algorithms are much less common than their discrete-time counterparts, and often only apply to algorithms with fairly simple dynamics \cite{Nevelson1976,Sen1978,Chen1982,Yin1993,Chen1994}. 
This being said, the continuous-time setting has recently been the subject of renewed attention. In particular, \cite{Sirignano2017a} introduced `{stochastic gradient descent in continuous time}' for recursive estimation of the parameters of a fully observed diffusion process, establishing a.s. convergence and a central limit theorem \cite{Sirignano2020a}. This algorithm has since been generalised in various directions \cite{Surace2019,Bhudisaksang2021,Sharrock2021a,Sharrock2020,Sharrock2020a}.

In this article, we provide the asymptotic convergence rate of a continuous-time TTSA algorithm designed for stochastic bilevel optimisation problems in continuous-time. This algorithm represents the continuous-time counterpart of \cite{Hong2020}, while the main result represents a rigorous extension of \cite{Mokkadem2006} to the continuous-time setting. 
We also demonstrate how this algorithm can be applied to several continuous-time bilevel optimisation problems

 

\section{Two-Timescale Stochastic Approximation in Continuous Time} 
\label{sec:algorithm}

\subsection{The Algorithm}
We are interested in obtaining solutions $x^{*}\in\mathbb{R}^{d_1}$ and $y^{*}\in\mathbb{R}^{d_2}$ which satisfy the first order stationary conditions corresponding to the bilevel optimisation problem in \eqref{eq:bilevel}, viz $\nabla \Phi(x^{*}) = 0$, $\nabla_{y}g(x^{*},y^{*}) = 0$, 
based on noisy observations of the gradients of the two objective functions. Perhaps the simplest way to tackle this problem is to simultaneously update $(x_t)_{t\geq 0}$ based on stochastic estimates of $\nabla_{x} \Phi(x_t)_{t\geq 0} \approx \nabla_{x} f(x_t,y_t)_{t\geq 0}$, and $(y_t)_{t\geq 0}$ based on stochastic estimates of $\nabla_{y}g(x_t,y_t)_{t\geq 0}$. However, while this approach has provable convergence guarantees and may perform well in practice \cite{Yang2019,Sharrock2020a}, it naturally ignores the dependence of the outer-level objective function $\Phi(x) = f(x,y^{*}(x))$ on $x$ in second argument. Thus, following \cite{Ghadimi2018}, 
we will instead consider an algorithm based on a more refined approximation of the gradient of the outer-level objective function, given by
\begin{equation}
\bar{\nabla}_{x} f(x,y) = \nabla_{x} f(x,y) - \nabla_{xy}^2 g(x,y) \left[ \nabla_{yy}^2 g(x,y)\right]^{-1} \nabla_{y} f(x,y).
\end{equation}
We will assume that our noisy gradient observations originate from the stochastic differential equations (SDEs)
\begin{align}
\mathrm{d}h_t^{(1)} &= \bar{\nabla}_{x} f(x_t,y_t)\mathrm{d}t +  \eta^{(1)}_t(x_t,y_t) \mathrm{d}t + \sigma_t^{(1)}(x_t,y_t) \mathrm{d}w_t^{(1)} \label{eq:obs1} \\
\mathrm{d}h_t^{(2)} &= \nabla_{y} g(x_t,y_t)\mathrm{d}t  + \eta^{(2)}_t(x_t,y_t) \mathrm{d}t +  \sigma_t^{(2)}(x_t,y_t) \mathrm{d}w_t^{(2)} \label{eq:obs2}
\end{align}
where, for $i=1,2$, $\eta_t^{(i)}:\mathbb{R}^{d_1}\times\mathbb{R}^{d_2}\rightarrow\mathbb{R}^{d_i}$, $\sigma_t^{(i)}:\mathbb{R}^{d_1}\times\mathbb{R}^{d_2}\rightarrow\mathbb{R}^{d_i\times d_i}$ are Borel measurable functions, and $w_t^{(i)}$ are standard $\mathbb{R}^{d_i}$-valued Brownian motions on a complete probability space $(\Omega,\mathcal{F},\mathbb{P})$. 
This naturally leads to the following algorithm for solving the bilevel optimisation problem \eqref{eq:bilevel}: recursively update $(x_t)_{t\geq 0}$ and $(y_t)_{t\geq 0}$ according to
\begin{align}
\mathrm{d} x_t & = \gamma_t^{(1)} [ \bar{\nabla}_{x} f(x_t,y_t)\mathrm{d}t + \eta^{(1)}_t(x_t,y_t) \mathrm{d}t + \sigma_t^{(1)}(x_t,y_t) \mathrm{d}w_t^{(1)}] \label{eq:SGD1} \\
\mathrm{d} y_t & = \gamma_t^{(2)} [ \nabla_{y} g(x_t,y_t)\mathrm{d}t  + \eta^{(2)}_t(x_t,y_t) \mathrm{d}t + \sigma_t^{(2)}(x_t,y_t) \mathrm{d}w_t^{(2)}] \label{eq:SGD2}
\end{align}
where $\gamma_t^{(i)}:\mathbb{R}_{+} \rightarrow\mathbb{R}^{d_i}_{+}$, $i=1,2$, are positive, non-increasing functions known as the {learning rates}. This represents the continuous-time analogue of the algorithm in \cite{Hong2020}. We will refer to it as TTSA in continuous time.

\subsection{Applications} \label{subsec:applications}
Many interesting problems can be seen as special cases of the stochastic bilevel optimisation problem in \eqref{eq:bilevel}, with noisy gradient estimates observed continuously in time as in \eqref{eq:obs1} - \eqref{eq:obs2}. We discuss several such problems below. 

\subsubsection{Model Agnostic Meta-Learning in Continuous-Time Models}
An important problem in machine learning is to obtain models which can quickly adapt to new tasks, based on limited new training data, or \emph{meta-learning}. 
One popular approach for performing meta-learning is model-agnostic meta-learning (MAML) \cite{Finn2017}. 
We are interested in MAML in the context of statistical estimation for continuous time models; see also \cite{Arango2021}. In particular, let us consider a collection of $\mathbb{R}^{n}$ valued diffusion processes $(z_t^{i})_{t\geq 0}$ defined by
\begin{equation}
\mathrm{d}z_t^{i} = h^{(i)}(z_t^{i})\mathrm{d}t + \varepsilon \mathrm{d} v_t^{i}, \label{eq:SDE}
\end{equation}
where $h^{i}:\mathbb{R}^{n}\rightarrow\mathbb{R}^{n}$, $\varepsilon \in\mathbb{R}^{n\times n}$, and $(v_t^{i})_{t\geq 0}$ are $\mathbb{R}^{n}$-valued standard Brownian motions. We will suppose that $(z_t^{i})_{t\geq 0}$ each possess a unique invariant measure, denoted $\pi^{i}(\mathrm{d}z)$. Given data $(z_t^{i})_{t\geq 0}^{i=1:N}$, we are then interested in learning a model for the drift function $h(z,\theta)$, parametrised by $\theta\in\mathbb{R}^p$, which quickly adapts to new observations. Following \cite{Rajeswaran2019,Hong2020}, we can formulate this as a bilevel optimisation problem as in \eqref{eq:bilevel}, namely, 
\begin{equation}
\min_{\theta\in\mathbb{R}^p} \Phi(\theta) = \sum_{i=1}^N \mathcal{L}^{(i)}(\theta^{*(i)}(\theta))~~~\text{subject to}~~~ \theta^{*(i)}(\theta)\in\argmin_{\theta^{(i)}\in\mathbb{R}^p} \left[ \mathcal{J}^{(i)}(\theta,\theta^{(i)}) \right],
\end{equation}
where $\theta$ is the shared model parameter, $\theta^{(i)}$ is the adaptation of $\theta$ to $(z_t^{i})_{t\geq 0}$, $\smash{\mathcal{J}^{(i)}(\theta,\theta^{(i)}) = \mathcal{L}^{(i)}(\theta^{(i)})+ \frac{\lambda}{2}||\theta^{(i)} - \theta||^2}$, and $\mathcal{L}^{(i)}$ is the asymptotic log-likelihood of $(z_t^{i})_{t\geq 0}$. In general, we cannot compute $\bar{\nabla}_{\theta} \mathcal{L}^{(i)}(\theta)$ or $\nabla_{\theta^{(i)}} \mathcal{J}^{(i)}(\theta,\theta^{(i)})$. We can, however, obtain stochastic estimates, in a similar form to \eqref{eq:obs1} - \eqref{eq:obs2}. This problem is thus amenable to continuous-time TTSA.

\subsubsection{Optimal Continuous-Time Optimisation Algorithms}
One of the most ubiquitous problems in computational statistics and machine learning is that of global optimisation, that is, finding $z^{*} = \argmin_{z\in\mathbb{R}^n} V(z)$ for a possibly non-convex function $V:\mathbb{R}^{n}\rightarrow\mathbb{R}$. One approach to this problem is simulated annealing, which recursively defines a sequence of estimates $(z_t)_{t\geq 0}$ as the solution of the SDE
\begin{equation}
\smash{\mathrm{d}z_t = -\nabla V(z_t)\mathrm{d}t + \sqrt{2\beta_t}\mathrm{d}v_t} \label{eq:SA}
\end{equation}
where $(v_t)_{t\geq 0}$ is a $\mathbb{R}^{d}$-valued standard Brownian motion, and $\beta_{t}:\mathbb{R}_{+}\rightarrow\mathbb{R}^n_{+}$ is a positive, non-increasing function known as the annealing schedule. These dynamics have clear connection with the related problem of sampling from an unnormalised probability distribution on $\mathbb{R}^d$.
Indeed, if $\beta_t=\beta$ is constant, this SDE admits as its unique invariant law the Gibbs distribution with density $\pi_{\beta}(z) \propto \exp [-\beta^{-1}V(z)]$. 
Meanwhile, if $\beta_t$ converges to zero at an appropriate rate, then $\pi_{\beta_t}(x)$ converges weakly to the set of global minima based on Laplace's principle \cite{Geman1986,Chiang1987}. 

One can also consider other 
versions of \eqref{eq:SA} with improved convergence properties.
 One such diffusion is obtained via the addition of a perturbation $\gamma$, satisfying $\nabla \cdot (\pi_{\beta} \gamma) = 0$, into the drift of \eqref{eq:SA}.
It is well known that this non-reversible diffusion converges faster to $\pi_{\beta}$, and admits a lower asymptotic variance, than its reversible counterpart \cite{Hwang2005,Duncan2016}. 
One can also consider a system of diffusions, interacting via a matrix $A$, which also exhibit improved convergence properties \cite{Borovykh2021a}. Given these results, it is natural to ask whether is it possible to determine an optimal perturbation, or an optimal interaction matrix.
One can view this as a simplified version of the bilevel optimisation problem in \eqref{eq:bilevel}. Suppose we parametrise $\gamma(z):=\gamma(\theta,z)$ or $A(z) := A(\theta,z)$, for some $\theta\in\mathbb{R}^p$. We then seek
\begin{equation}
{\theta^{*} = \argmin_{\theta\in\mathbb{R}^p}\Phi(\theta) = U(\theta,z^{*})~~~\text{subject to}~~~z^{*}\in\argmin_{z\in\mathbb{R}^n} V(z)} \label{eq:bilevel3}
\end{equation}
where $U$ encodes the asymptotic variance, the rate of convergence, or any other relevant constraints (e.g., the cost of interactions). One can obtain noisy estimates for the gradients of the two objective functions, in a slightly more general form of \eqref{eq:obs1} - \eqref{eq:obs2}, and thus solve \eqref{eq:bilevel3} using a slight modification of the continuous-time TTSA algorithm in \eqref{eq:SGD1} - \eqref{eq:SGD2}. We note that here the inner-level problem is independent of the outer-level problem, leading to the simplification $\bar{\nabla}_{\theta}U(\theta,z) = \nabla_{\theta}U(\theta,z)$.

\subsubsection{Online Parameter Estimation and Optimal Sensor Placement in Partially Observed Diffusions}
Partially observed diffusion processes are frequently used to model stochastic dynamics in engineering, economics, and the natural sciences. Such processes are described by two SDEs: a signal equation for the latent $\mathbb{R}^{n_x}$-valued signal process $(x_t)_{t\geq 0}$, and an observation equation for the $\mathbb{R}^{n_y}$-valued observation process $(y_t)_{t\geq 0}$ \cite{Bain2009}. Typically, the signal equation will depend on $\theta\in\mathbb{R}^{n_p}$, a vector of model parameters which must be estimated from data, preferably in an online fashion \cite{Surace2019}. In addition, the observation equation may depend on $\boldsymbol{o} = (\boldsymbol{o}_i)_{i=1}^{n_y}$, a set of $n_y$ measurement sensors in $\mathbb{R}^{n_{\boldsymbol{o}}}$, which can often be moved to improve the optimal estimate of the latent state \cite{Athans1972}.

It is natural to ask whether one can jointly perform online parameter estimation and optimal sensor placement. This problem can be precisely formulated as a stochastic bilevel optimisation problem, viz
\begin{equation}
{\min_{\theta\in\mathbb{R}^{n_{\theta}}} \Phi(\theta) = \mathcal{L}(\theta,\boldsymbol{o}^{*}(\theta))~~~\text{subject to}~~~\boldsymbol{o}^{*}(\theta) \in \argmin_{\smash{\boldsymbol{o}\in(\mathbb{R}^{n_{\boldsymbol{o}}})^{n_y}}} \mathcal{J}(\theta,\boldsymbol{o}),}
\end{equation} 
where $\mathcal{L}(\theta,\boldsymbol{o})$ is the asymptotic log-likelihood of the observations, and $\mathcal{J}(\theta,\boldsymbol{o})$ the the asymptotic covariance of the optimal state estimate. As shown rigorously in \cite{Sharrock2020a}, one can thus solve this problem using a continuous-time TTSA algorithm, using a generalised version of \eqref{eq:SGD1} - \eqref{eq:SGD2},


\section{Main Results}
\label{sec:main}
In this section, we present results on the a.s. convergence and convergence rate of the continuous-time TTSA algorithm.

\subsection{Almost Sure Convergence}
\label{subsec:assumptions}
\begin{assumption} \label{assumption:learningrate}
The learning rates are given by $\smash{\gamma_t^{(1)} = \gamma_0^{(1)} (\delta_1 + t)^{-\eta_1}}$ and $\smash{\gamma_t^{(2)} = \gamma_0^{(2)} (\delta_2 + t)^{-\eta_2}}$,  
where $\gamma_0^{(1)},\gamma_0^{(2)}, \delta_1,\delta_2\in\mathbb{R}_{+}$, and $\eta_1,\eta_2\in(\frac{1}{2},1)$ satisfy $\eta_2<\eta_1$. 
\end{assumption}

\begin{assumption} \label{assumption:outer}
The functions $\Phi:\mathbb{R}^{d_1}\rightarrow \mathbb{R}$ and $f:\mathbb{R}^{d_1}\times \mathbb{R}^{d_2}\rightarrow\mathbb{R}$ satisfy: (i) for all $x\in\mathbb{R}^{d_1}$, $\nabla_{x}f(x,\cdot)$ and $\nabla_{y}f(x,\cdot)$ are Lipschitz continuous w.r.t. $y$, (ii) for all $y\in\mathbb{R}^{d_2}$, $\nabla_{y} f(\cdot,y)$ is Lipschitz continuous w.r.t. $x$, (iii) for all $x\in\mathbb{R}^{d_1}$, $y\in\mathbb{R}^{d_2}$, $||\nabla_{y} f(x,y)||$ is bounded above. 
\end{assumption}

\begin{assumption} \label{assumption:inner}
The function $g:\mathbb{R}^{d_1}\times\mathbb{R}^{d_2}\rightarrow \mathbb{R}$ satisfies: (i) for all $x\in\mathbb{R}^{d_1}$, $y\in\mathbb{R}^{d_2}$, $g(x,y)$ is twice continuously differentiable in $(x,y)$, (ii) for all $x\in\mathbb{R}^{d_1}$, $\nabla_{y} g(x,\cdot)$ is Lipschitz continuous w.r.t. $y$, (iii) for all $x\in \mathbb{R}^{d_1}$, $y\in\mathbb{R}^{d_2}$, $\nabla_{yy}^2 g(x,y)\succcurlyeq\mu_{g}I$, (iv) for all $x\in\mathbb{R}^{d_1}$, $\nabla_{xy}^2g(x,\cdot)$ and $\nabla_{yy}^2g(x,\cdot)$ are Lipschitz continuous w.r.t. $y$, (v) for all $y\in\mathbb{R}^{d_2}$, $\nabla_{xy}^2 g(\cdot,y)$ and $\nabla_{yy}^2g(\cdot,y)$ are Lipschitz continuous w.r.t $x$, (vi) for all $x\in\mathbb{R}^{d_1}$, ${y}\in\mathbb{R}^{d_2}$, $||\nabla_{xy}^2 g(x,y)||$ is bounded above. 
\end{assumption}
\vspace{.5mm}
\begin{assumption} \label{assumption:bound}
The algorithm iterates $(x_t)_{t\geq 0}$ and $(y_t)_{t\geq 0}$ are almost surely bounded.
\end{assumption}

\vspace{.5mm}
\begin{theorem} \label{theorem:convergence}
Suppose that Assumptions \ref{assumption:learningrate} - \ref{assumption:bound} are satisfied. Then, almost surely, $\lim_{t\rightarrow\infty} x_t = x^{*}$  and $\lim_{t\rightarrow\infty} y_t = y^{*}$.
\end{theorem}
\vspace{-4mm}
\begin{proof}
See \cite{Sharrock2020a}. 
\end{proof}

\subsection{Convergence Rate}
\begin{assumption} \label{assumption:outer2}
There exists a neighbourhood $\mathcal{U}_{(x^{*},y^{*})}$ of $(x^{*},y^{*})$ such that $\bar{\nabla}^2_{xx} f(x,y)\succcurlyeq\mu_{f}I$ for all $(x,y)\in\mathcal{U}_{(x^{*},y^{*})}$, where $\bar{\nabla}_{xx}^2f(x,y) = \nabla_{x} [\bar{\nabla}_{x}f(x,y)] -  \nabla_{xy}^2g(x,y) [\nabla_{yy}^2 g(x,y)]^{-1} \nabla_{y} [\bar{\nabla}_{x}f(x,y)]$.
\end{assumption}

\begin{assumption} \label{assumption:noise1}
The functions $\eta_t^{(i)}:\mathbb{R}^{d_1}\times\mathbb{R}^{d_2}\rightarrow\mathbb{R}^{d_i}$ satisfy $\eta_t^{(i)}(x_t,y_t) = o((\gamma_t^{(1)})^{\frac{1}{2}})$.
\end{assumption}

\begin{assumption}  \label{assumption:noise2}
The functions $\sigma_t^{(i)}:\mathbb{R}^{d_1}\times\mathbb{R}^{d_2}\rightarrow\mathbb{R}^{d_i\times d_i}$ satisfy $\smash{\lim_{t\rightarrow\infty} \sigma^{(1)}_t[\sigma^{(1)}_t]^T= \Gamma_{11}}$, $\smash{\lim_{t\rightarrow\infty} \sigma^{(2)}_t[\sigma^{(2)}_t]^T= \Gamma_{22}}$, and $\smash{\lim_{t\rightarrow\infty} \sigma^{(1)}_t[\sigma^{(2)}_t]^T= \Gamma_{12}}$. 
\end{assumption}

\vspace{.5mm}
\begin{theorem} \label{theorem:convergence_rate}
Suppose that Assumptions \ref{assumption:learningrate} - \ref{assumption:bound} and \ref{assumption:outer2} - \ref{assumption:noise2} are satisfied. Define $\Gamma^{*}_{11} = \Gamma_{11}(x^{*},y^{*})$, $\Gamma_{22}^{*}=\Gamma_{22}(x^{*},y^{*})$, and $\Gamma^{*}_{12} = \Gamma_{12}(x^{*},y^{*})$. In addition, define $A_{11} = -\nabla_{x}\bar{\nabla}_{x}f(x^{*},y^{*})$, $A_{12} = -\nabla_{y}\bar{\nabla}_{x}f(x^{*},y^{*})$, $A_{21} = -\nabla_{x}\nabla_{y}g(x^{*},y^{*})$, $A_{22} = -\nabla_{y}\nabla_{y}g(x^{*},y^{*})$, and $H = A_{11} - A_{12}A_{22}^{-1}A_{21}$. Finally, define 
\begin{align}
\Sigma_{x}& = \int_0^{\infty} \exp\left( Ht\right) \left(\Gamma_{11}^{*} + A_{12}A_{22}^{-1}\Gamma_{22}^{*} [A_{22}^{-1}]^T A_{12}^T - \Gamma_{12}^{*}[A_{22}^{-1}]^T A_{12}^T - A_{12}A_{22}^{-1} \Gamma_{21}^{*}\right) \exp\left( Ht\right) \mathrm{d}t \\
\Sigma_{y}& = \int_0^{\infty} \exp\left(A_{22}t\right) \Gamma_{22}^{*} \exp\left(A_{22}t\right)\mathrm{d}t.
\end{align}
Then we have that 
\begin{equation}
\begin{pmatrix}
(\gamma_t^{(1)})^{-\frac{1}{2}} ({x}_t-x^{*}) \\ (\gamma_t^{(2)})^{-\frac{1}{2}} (y_t - y^{*})
\end{pmatrix} 
\stackrel{\mathcal{D}}{\longrightarrow} \mathcal{N}\left(0,\begin{pmatrix} \Sigma_{x} & 0 \\ 0 & \Sigma_{y}\end{pmatrix} \right).
\end{equation}
\end{theorem}
\vspace{-5mm}
\begin{proof}
See \cite{Sharrock2022b}
\end{proof}

\section{Discussion} \label{sec:discussion}
There are many interesting directions in which the results presented in this article can be extended. Perhaps the most natural extension is to consider the case in which the noisy observations \eqref{eq:obs1} - \eqref{eq:obs2} depend on an additional ergodic continuous-time Markov process $(z_t)_{t\geq 0}$ with unique invariant law $\mu(\mathrm{d}z)$. That is, 
\begin{align}
\mathrm{d}h_t^{(1)} &= F(x_t,y_t,z_t) \mathrm{d}t + \eta^{(1)}_t(x_t,y_t) \mathrm{d}t + \sigma_t^{(1)}(x_t,y_t) \mathrm{d}w_t^{(1)}  \label{eq:obs1.1} \\
\mathrm{d}h_t^{(2)} &= G(x_t,y_t,z_t) \mathrm{d}t  + \eta^{(2)}_t(x_t,y_t) \mathrm{d}t + \sigma_t^{(2)}(x_t,y_t) \mathrm{d}w_t^{(2)} \label{eq:obs2.1}
\end{align}
where $\mathbb{E}_{\mu}[F(x,y,z)] = \bar{\nabla}_xf(x,y)$ and $\mathbb{E}_{\mu}[G(x,y,z)] = {\nabla}_y g(x,y)$. Indeed, in practice, this is the form of the noisy observations for the problems discussed in Section \ref{subsec:applications}. Under suitable conditions, one can establish a.s. convergence in this case, as in \cite{Sharrock2020a}[Section 2.2], by appealing to classical results on the solutions of a related Poisson equation \cite{Pardoux2001}. It is also possible to establish a central limit theorem, using Theorem \ref{theorem:convergence_rate} and the tools recently developed in \cite{Sirignano2020a}. The proof of this result will appear in a forthcoming paper \cite{Sharrock2022a}. More generally, it is of interest to extend other results on discrete-time, two-timescale stochastic approximation \cite{Dalal2018,Gupta2019,Kaledin2020,Doan2021a,Doan2022b} to the continuous-time setting. 

In practice, it is evident that any stochastic approximation scheme in continuous time must be discretised. Thus, when designing statistical learning algorithms for continuous-time (bilevel) optimisation problems, it is natural to ask why we would prefer to use a discrete-time approximation of continuous-time TTSA over the traditional approach, which first discretises the continuous-time model, and then applies discrete-time TTSA. In certain examples (e.g., value function estimation in continuous-time models), it is known that the continuous-time approach has certain advantages over the discrete-time approach \cite{Sirignano2017a}. We believe, however that further rigorous investigation into this question remains an important direction for future research.



\bibliography{two-timescale-refs}
\bibliographystyle{icml2022}

\end{document}